\documentclass[letterpaper, 10 pt, conference]{ieeeconf}  

\IEEEoverridecommandlockouts                              

\usepackage{amsmath, amsfonts}
\usepackage{color}
\usepackage{amssymb}
\usepackage{mathtools}
\usepackage{csquotes}
\usepackage{quotes}
\usepackage{balance}

\newtheorem{theorem}{Theorem}
\newtheorem{lemma}{Lemma}

\usepackage{rotating}
\usepackage{adjustbox}

\usepackage{calc}

\usepackage{longtable}
\usepackage{multirow}

\usepackage[disable]{todonotes} 

\usepackage{graphicx}

\usepackage{placeins}

\usepackage{tabularx}

\usepackage{lmodern,textcomp}

\usepackage{enumerate}

\usepackage{xspace}

\usepackage{algorithm}
\usepackage{algpseudocode}

\newcommand{\inv}{^{-1}\xspace}

\newcommand{\blambda}{\boldsymbol{\lambda}\xspace}

\newcommand{\bxi}{\boldsymbol{\xi}\xspace}
\newcommand{\bs}{\textbf{s}\xspace}
\renewcommand{\t}{^\textsf{T}\xspace}

\newcommand{\away}[1]{}

\newcommand{\R}{\mathbb{R}\xspace}

\newcommand{\N}{\mathbb{N}\xspace}

\newcommand{\cA}{\mathcal{A}\xspace}
\newcommand{\cB}{\mathcal{B}\xspace}

\newcommand{\cL}{\mathcal{L}\xspace}

\newcommand{\cJ}{\mathcal{J}\xspace}
\newcommand{\cU}{\mathcal{U}\xspace}

\newcommand{\cI}{\mathcal{I}\xspace}

\newcommand{\cO}{\mathcal{O}\xspace}

\newcommand{\bA}{\textbf{A}\xspace}
\newcommand{\bJ}{\textbf{J}\xspace}

\newcommand{\tbM}{\tilde{\textbf{M}}\xspace}
\newcommand{\tbH}{\tilde{\textbf{H}}\xspace}
\newcommand{\tbJ}{\tilde{\textbf{J}}\xspace}
\newcommand{\tbg}{\tilde{\textbf{g}}\xspace}

\newcommand{\bZ}{\textbf{Z}\xspace}

\newcommand{\bx}{\textbf{x}\xspace}

\newcommand{\bb}{\textbf{b}\xspace}

\newcommand{\bc}{\textbf{c}\xspace}
\newcommand{\bg}{\textbf{g}\xspace}

\newcommand{\tbc}{\tilde{\textbf{c}}\xspace}

\newcommand{\hbx}{\hat{\textbf{x}}\xspace}

\newcommand{\bz}{\textbf{z}\xspace}
\newcommand{\be}{\mathbf{1}\xspace}
\newcommand{\bei}[1]{{\textbf{e}}\xspace}

\newcommand{\bM}{\textbf{M}\xspace}
\newcommand{\bN}{\textbf{N}\xspace}

\newcommand{\bH}{\textbf{H}\xspace}
\newcommand{\bI}{\textbf{I}\xspace}

\newcommand{\bO}{\textbf{0}\xspace}

\newcommand{\tbx}{\tilde{\textbf{x}}\xspace}

\newcommand{\tol}{{\textsf{tol}}\xspace}

\newcommand{\commentout}[1]{}

\usepackage{datetime,url}

\usepackage{cite}

\title{\LARGE \bf
Solving Problems with Inconsistent Constraints with a Modified Augmented Lagrangian Method}

\author{Martin P. Neuenhofen$^{1}$ and Eric C. Kerrigan$^{2}$
\thanks{$^{1}$Martin P. Neuenhofen is with the Deparment of Electrical \& Electronic Engineering, Imperial College London, SW7 2AZ London, UK, e-mail: m.neuenhofen19@imperial.ac.uk
        {\tt\small www.MartinNeuenhofen.de}}%
\thanks{$^{2}$Eric C. Kerrigan is with the Deparment of Electrical \& Electronic Engineering and Department of Aeronautics, Imperial College London, SW7 2AZ London, UK, e-mail: e.kerrigan@imperial.ac.uk
        {\tt\small www.imperial.ac.uk/people/e.kerrigan}}%
}

\begin{document}

\sloppy
\maketitle
\thispagestyle{empty}
\pagestyle{empty}

\begin{abstract}

We present a numerical method for the minimization of constrained optimization problems where the objective is augmented with large quadratic penalties of inconsistent equality constraints. Such objectives arise from quadratic integral penalty methods for the direct transcription of optimal control problems.
The Augmented Lagrangian Method (ALM) has a number of advantages over the Quadratic Penalty Method (QPM).
However, if the equality constraints are inconsistent, then ALM might not converge to a point that minimizes the bias of the objective and penalty term. Therefore,  we present a modification of ALM that fits our purpose.
We prove convergence of the modified method and bound its local convergence rate by that of the unmodified method.
Numerical experiments demonstrate that the modified ALM can minimize certain quadratic penalty-augmented functions faster than QPM, whereas the unmodified ALM converges to a minimizer of a significantly different problem.
\end{abstract}

\newcommand{\abbTab}{Table}
\newcommand{\abbtab}{Table}
\newcommand{\abbAlg}{Algorithm}
\newcommand{\abbalg}{Algorithm}
\newcommand{\abbthm}{Theorem}
\newcommand{\abbThm}{Theorem}
\newcommand{\abbsec}{Section}
\newcommand{\abbSec}{Section}

\newcommand{\pval}{\omega\xspace} 	
\newcommand{\pmval}{\rho\xspace} 	
\newcommand{\wquad}{\alpha\xspace} 	
\newcommand{\tquad}{\tau\xspace} 	

\newcommand{\bdual}{\blambda} 						
\newcommand{\hbA}{\hat{\boldsymbol{\bA}}\xspace} 	
\newcommand{\bbar}{\boldsymbol{\eta}\xspace} 		

\newcommand{\eqnKKT}{(KKT)\xspace}
\newcommand{\eqnKKTa}{(KKT1)\xspace}
\newcommand{\eqnKKTb}{(KKT2)\xspace}

\section{Introduction}

\subsection{Problem Statement}
This paper describes and analyzes a modified augmented Lagrangian method (MALM) for the numerical solution of a \textit{quadratic penalty program}:
\begin{align*}
\min_{\bx \in \cB} \ \Phi_{\pval}(\bx)&:= f(\bx) + \frac{1}{2  \pval}  \|c(\bx)\|_2^2\,, \tag{QPP}\label{eqn:CQPP}\\
\cB&:= \lbrace \bx \in \R^n\, \vert\, g(\bx)\geq \bO \rbrace\,,
\end{align*}
where $f : \R^n \rightarrow \R$, $c : \R^n \rightarrow \R^m$, $g : \R^n \rightarrow \R^p$ are possibly non-convex and nonlinear functions; $\geq$ is meant for each vector component; $\cB$ is the feasible set; $m,n,p \in \N$ are dimensions; $\pval \in \R_{>0}$ is part of the problem data.

We define the associated \emph{Lagrangian function} $\cL(\bx,\bdual,\bbar):=f(\bx)-\bdual\t\cdot c(\bx) - \bbar\t\cdot g(\bx)$, with Lagrange multipliers $\bdual \in \R^m$, $\bbar \in \R^p_{\geq 0}$. We explain the meaning of a Lagrange multiplier $\bdual$ below.

\subsubsection{Relation to Constrained Programs (CP)}\label{sec:Intro:Relate}
When $\pval>0$ is close to zero then the penalty forces $c(\bx)\approx \bO$, provided such a point exists. Hence, the problem may be considered to be related to:
\begin{equation*}
\min_{\bx \in \cB} \qquad f(\bx) \quad\text{s.t.  }c(\bx)=\bO
\tag{CP}\label{eqn:CP}
\end{equation*}
In~\eqref{eqn:CP}, $c,g$ are equality and inequality constraint functions with Lagrange multipliers $\bdual \in \R^m$, $\bbar \in \R^p_{\geq 0}$.

\subsubsection{Inconsistency}
\eqref{eqn:CP} only makes sense when $c(\bx)=\bO$ is consistent. However, for the scope of this work we are particularly interested in the case when $c$ is inconsistent. Experiments show that for inconsistent $c$ the solution of~\eqref{eqn:CQPP} depends significantly on the value of $\pval$; cf. Section~\ref{sec:Intro:Motiv} and Figure~\ref{fig:ocp_study}.

\subsubsection{Optimality Conditions}
From \cite[Thm~12.1]{NumOpt}:
\begin{align*}
	&\underbrace{\nabla f(\bx) - \nabla c(\bx) \cdot \frac{-1}{\pval} \cdot c(\bx)}_{\equiv \nabla \Phi_{\pval}(\bx)} - \nabla g(\bx) \cdot \bbar = \bO\tag{KKT1'}\label{eqn:KKT1'}\\
	&g_i(\bx)=0 	 \text{ and }\bbar_i\geq 0 		\qquad \forall i \in \cA\tag{KKT2a}\\
	&g_i(\bx)>0  	 \text{ and }\bbar_i =   0     	\qquad \forall i \notin \cA\tag{KKT2b}
\end{align*}
where $\cA \subseteq \lbrace 1,\dots,p\rbrace$ is the \textit{active set}, and $g_i$ is the $i^\text{th}$ component of the vector $g(\bx)$.

Substituting $\bdual = \frac{-1}{\pval} \cdot c(\bx)$, we can re-express~\eqref{eqn:KKT1'}:
\begin{align}
	\nabla_\bx \cL(\bx,\bdual,\bbar)=\bO\,,\qquad c(\bx) + \pval \cdot \bdual =\bO
	\tag{KKT1}\label{eqn:KKT1}
\end{align}
\eqnKKT (i.e., \eqnKKTa and \eqnKKTb) determines $\bx,\bdual,\bbar$. \eqnKKT are the optimality conditions of~\eqref{eqn:CQPP} when $\pval>0$ and the optimality conditions of~\eqref{eqn:CP} when $\pval=0$.

\subsection{Motivation}

\subsubsection{Necessity of Tailored Solvers for~\eqref{eqn:CQPP}}\label{sec:Intro:Motiv:Num}
Minimizing~\eqref{eqn:CQPP} directly appears natural but, unless $c$ is affine, will result in many iterations. This is caused by the bad scaling of the penalties.

As a demonstration, consider the instance
\begin{subequations}
	\begin{align}
	f(\bx)&:= -x_1\label{eqn:CircleF}\,,\quad g(\bx) :=\begin{bmatrix} x_1\\
	x_2-x_1
	\end{bmatrix} \in \R^{2}\\
	c(\bx)&:= \begin{bmatrix}
	(x_1 +\varepsilon)^2 + x_2^2 - 2\\
	(x_1 -\varepsilon)^2 + x_2^2 - 2
	\end{bmatrix} \in \R^{2}\label{eqn:CircleC}
	\end{align}
	\label{eqn:Circle}%
\end{subequations}
with primal and dual initial guesses $\bx_0 :=[2 \quad 1]\t$ and $\bdual_0 :=\bO$, for $\varepsilon=0$. We discuss later with \abbtab~\ref{tab:Exp_Circ_iter} that minimization of~\eqref{eqn:CQPP} of~\eqref{eqn:Circle} with a direct minimization method takes $134$ iterations when $\pval=10^{-6}$. This is inefficient when compared to our later proposed MALM, which solves the same instance in only $39$ iterations.

\subsubsection{Relevant Instances of~\eqref{eqn:CQPP}}\label{sec:Intro:Motiv}
Integral penalty methods \cite{Balakrishnan68,Hager90,CDC2020_PBF} are an alternative to collocation methods for solving dynamic optimization problems. Integral penalty methods can solve dynamic optimization problems with singular arcs and high-index differential-algebraic path-constraints as a problem of form~\eqref{eqn:CQPP}. Consider the bang-singular example
\begin{equation*}
\begin{aligned}
&\min_{y,u} &\quad &J:=\int_0^5 \left(\,y(t)^2 + t\,u(t)\,\right)\,\mathrm{d}t,\\
&\text{s.t.} & y(0)&=0.5,\quad\dot{y}(t)=\frac{1}{2}y(t)^2+u(t)\,,\\
&  & y(t),u(t)& \in [-1,1]\quad\forall t \in [0,5]\,.
\end{aligned}\tag{OCP}\label{eqn:ExampleOCP}
\end{equation*}
In integral-penalty-methods, the idea is to force $y(t)^2/2+u-\dot{y}=0$ not only at collocation points\todo{2f}, but instead add an integral penalty $r=\int_0^5 \|y^2/2+u-\dot{y}\|_2^2\,\mathrm{d}t$ to the objective.

\begin{figure}
	\centering
	\includegraphics[width=0.75\linewidth]{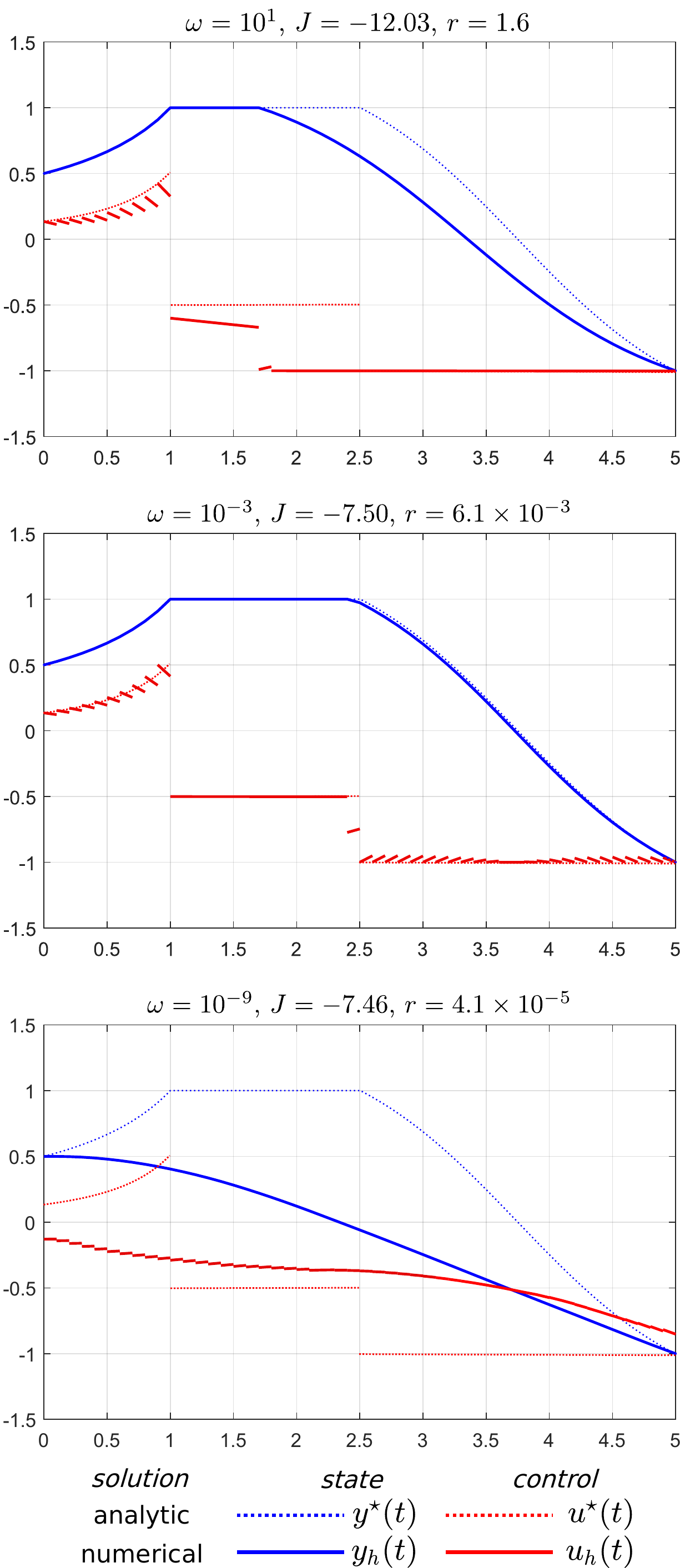}
	\caption{Numerical solution to~\eqref{eqn:ExampleOCP} for $N=40$ and different values of $\pval$.}
	\label{fig:ocp_study}
\end{figure}

Consider using continuous piecewise linear finite elements $y_h$ for $y$ and discontinuous ones $u_h$ for $u$ on a uniform mesh of $N \in \N$ intervals (mesh size $h=5/N$); represented with $\bx := [y_h(h),\dots, y_h(Nh), u^+_h(0), u^-_h(h),u^+_h(h) \dots, u^-_h(Nh)]\t \in \R^n$, $n:={3 N}$. $y_h(0)=0.5$ is fixed and removed from $\bx$. We can minimize a quadrature approximation of $J + \frac{1}{2\pval}r$ by solving an instance of~\eqref{eqn:CQPP}, where
\begin{subequations}
	\label{eqn:OCPdisc}
	\begin{align}
	f(\bx)&:= \sum_{j=1}^{Nq} \wquad_{j} \left( y_h(\tau_{j})^2 + \tau_{j} \, u_h(\tau_j) \right)\label{eqn:OCPdiscF}\\
	c(\bx)&:= \begin{bmatrix} \vdots \\
	\sqrt{\wquad_{j}} \big(y_h(\tau_j)^2/2 + u_h(\tau_j)-\dot{y}_h(\tau_{j}) \big)\\
	\vdots
	\end{bmatrix} \in \R^{m}\label{eqn:OCPdiscC}\\
	g(\bx)&:= \begin{bmatrix}
	\be-\bx\\
	\be+\bx
	\end{bmatrix} \in \R^{p}\label{eqn:OCPdiscG}
	\end{align}
\end{subequations}
with $q$ quadrature points $\tau_j$ and weights $\alpha_j>0$ per mesh-interval. $\pval \in \R_{>0}$ is ideally chosen in $\cO(1/N)$ \cite{Hager90,PBF}. Figure~\ref{fig:ocp_study} plots numerical solution $y_h,u_h$ for different values of $\pval$ against the analytic solution. The numerical solutions are vastly different for different $\pval$. Problem~\eqref{eqn:CP} is infeasible for~\eqref{eqn:OCPdisc}\todo{2j} because $c(\bx)\neq \bO$ $\forall \bx \in \cB$.

In conclusion: Problems~\eqref{eqn:CQPP} and~\eqref{eqn:CP} have different solutions. Solutions of~\eqref{eqn:CQPP} depend on $\pval$. Due to page limitations,\todo{2e} for a discussion on integral penalty methods, implementation, and choice of $\pval$, we refer to \cite{Hager90,Balakrishnan68,CDC2020_PBF}. \todo{3a}The experiments in \cite{PBF,CDC2020_PBF} present singular-arc and differential-algebraic optimal control problems where collocation methods fail to converge, but integral penalty methods converge.

\subsection{Literature Review}
We saw in Section~\ref{sec:Intro:Motiv:Num} that straightforward numerical minimization of~\eqref{eqn:CQPP} is inefficient due to bad scaling when $\pval$ is close to zero, hence necessitating tailored algorithms.

\subsubsection{Quadratic Penalty Method (QPM)}
QPMs compensate for the bad scaling by iteratively minimizing a sequence of problems~\eqref{eqn:CQPP}. Therein, $\pval$ is replaced by\todo{2i} a sequence of values $\lbrace\pmval_k\rbrace_{k \in \N_0}$ that converges to $\pval$ from above. Due to page limitations, \todo{2n}we refer to \cite{Courant43,SUMT} for details. Actually, these methods have been proposed for problem~\eqref{eqn:CP}, i.e.\ when $\pval=0$; but they can also be used for~\eqref{eqn:CQPP}. This is so because QPMs solve penalty problems of form~\eqref{eqn:CQPP}. QPMs can converge slowly due to bad scaling \cite{Murray71}.

\subsubsection{Augmented Lagrangian Method (ALM)}
ALMs have been developed as a replacement for QPMs when solving~\eqref{eqn:CP}. They work like QPMs, but augment $\Phi_{\pmval_k}(\bx)$ with the term $-\bdual_{k}\t \cdot c(\bx)$. This term with $\bdual_k \in \R^m$ creates the right amount of slope such that the inequality constrained minimizer of $\Psi_k$ eventually matches with the minimizer of~\eqref{eqn:CP}. Due to page limitations, \todo{2o}we refer to \cite[Alg.~17.3]{NumOpt} and the references therein for all details on how $\bdual_k \in \R^m$ is iteratively refined to achieve this. Convergence of $\bdual_k$ is asserted under suitable conditions \cite{Bertsekas1,Conn95}.

\subsubsection{Extensions of ALM to Inequality Constraints}
Originally, ALM treated only equality constraints \cite{Hestenes1,Powell1} by means of quadratic penalties of $c$ and update schemes for $\bdual$. \todo{1a}In this case, inequalities can be \emph{subjected}~\cite{Conn95}, i.e.\ minimize the sequence $\Psi_{k}$ subject to $\bx \in \cB$. Alternatively, penalty or barrier terms of $g$ can be \emph{augmented} \cite{Rockafellar1} with according update strategies for $\bbar$. Subjections are considered more efficient in practice than augmentations~\cite{Lancelot}. Augmentations can suffer from non-smooth, non-differentiable, or low-order smooth penalties/barriers, and can converge slower or less reliably.\todo{1b}

\subsubsection{Extensions of ALM to~\eqref{eqn:CQPP}}
Originally, ALM treated only~\eqref{eqn:CP} as opposed to~\eqref{eqn:CQPP}. The work \cite{SHARIFF2003257} proposes a modified scheme (MALM) for~\eqref{eqn:CQPP} when $f$ is quadratic, $c$ linear, $\Phi_{\pval}$ convex, and $\cB=\R^n$. They prove global convergence of the their scheme. Our previous work in \cite{CDC2020_MALM} extended MALM to problems where $f$ is general, $c$ nonlinear, without additional inequality constraints $g$ \todo{1d}. Also, there is no convergence proof yet in the literature for the case when $f,c$ are general, regardless of the presence of $g$. This paper will present such proofs \todo{1e}.

\subsection{Challenges}
Our goal is in devising a method that solves~\eqref{eqn:CQPP} by solving a sequence of penalty problems with moderate penalty parameter $\pmval \gg \pval$, and prove its convergence. In the limit $\pval \rightarrow 0$, MALM should match ALM due to the relation of the problems~\eqref{eqn:CQPP} and~\eqref{eqn:CP} as described in Section~\ref{sec:Intro:Relate}.

Proving convergence for non-convex $\Phi_{\pval}$ is challenging because solutions of sub-problems may be non-unique and hence alternating. Convergence of $\bbar_k$ may be challenging to prove because the solution $\bbar$ of \eqnKKT may be non-unique. We will assert uniqueness of $\bbar$ from a strict complementarity assumption. Striking the right balance between mild assumptions and strong convergence assertions appears non-trivial in this context.

\subsection{Contributions}
We present MALM for general functions $f,c,g$ (Algorithm~1). We prove convergence for the case when $f,c$ are twice continuously differentiable and $g$ is linear (Theorem~\ref{thm:globconv}). Furthermore, we give a local rate-of-convergence result for the case when $f,c,g$ are twice local Lipschitz-differentiable (Theorem~\ref{thm:locconv}).

Theorem~1 is not easily extendable to nonlinear $g$ because it uses a result for ALM on~\eqref{eqn:CP} for linear~$g$. Theorem~2 works for general $g$ but assumes convergence and Lipschitz-continuous second derivatives of $f,c,g$. In the iteration limit, convergence of ALM can only be guaranteed to be at least at a\todo{1c} linear rate \cite[Thm~17.6]{NumOpt}. Likewise, our rate-of-convergence result for MALM asserts only a linear rate. However, this linear rate is slightly better than the linear rate of ALM. Hence, our work draws connections between the rate of convergence between MALM and ALM.

\subsection{Structure of the Paper}
Section~II derives the proposed algorithm. Section~III presents the convergence analysis. Section~IV gives numerical experiments.

\section{Derivation of the Algorithm of MALM}
\label{sec:MALM}
MALM is a solution method for~\eqref{eqn:CQPP}. MALM has been presented\todo{2a} in \cite{SHARIFF2003257} for the special case when $f$ is quadratic, $c$ linear, and $\Phi_{\pval}$ convex. The method has been presented for the case where $f,c$ are general in \cite{CDC2020_MALM} but without inequality constraints and without a convergence analysis\todo{1f}. Here, we derive MALM for general nonlinear non-convex $f,c,g$, and in a stronger relation\todo{2b} to its origins in ALM \cite{Hestenes1,Powell1}. For the method presented here, we give global and local convergence proofs.

The derivation poses an auxiliary problem, applies ALM to it, and then eliminates variables.

\subsection{Auxiliary Problem}
The following problem is equivalent to~\eqref{eqn:CQPP} but of the form~\eqref{eqn:CP}:
\begin{subequations}
	\begin{align}
		&\operatornamewithlimits{min}_{\hbx:=(\bx,\bxi) \in \R^{(n+m)}} 	&\quad 	\hat{f}(\hbx)&:=f(\bx) + \frac{\pval}{2}  \|\bxi\|_2^2 \label{eqn:SubstObjective}	\\
		&\text{s.t.} 	&		\hat{c}(\hbx)&:=c(\bx) + \pval  \bxi =\bO\,, \label{eqn:SubstConstraints}\\[3pt]
		& 				 	&		\hat{g}(\hbx)&:=g(\bx)\geq \bO\,.
	\end{align}
	\label{eqn:Subst}%
\end{subequations}
The optimality conditions of~\eqref{eqn:Subst} are \eqnKKTb and
\begin{subequations}
	\label{eqn:KKT1_Subst}
	\begin{align}
	\begin{bmatrix}
	\nabla f(\bx)\\
	\pval\bxi
	\end{bmatrix} - \begin{bmatrix}
	\nabla c(\bx)\\
	\pval  \bI
	\end{bmatrix}  \bdual -\begin{bmatrix}
	\nabla g(\bx)\\
	\bO
	\end{bmatrix} \bbar &=\bO\\
	c(\bx) + \pval  \bxi &= \bO\,. \label{eqn:5b}
	\end{align}
\end{subequations}

\subsection{Augmented Optimality System}
Since~\eqref{eqn:Subst} is of form~\eqref{eqn:CP}, we can apply ALM with augmented inequality constraints\todo{2c} as in \cite{Conn95}. To this end, we introduce an auxiliary vector $\bz \in \R^m$ and a moderate penalty parameter $\pmval>0$. These are added to~\eqref{eqn:SubstConstraints} and in the gradient of the Lagrangian function:
\begin{subequations}
	\begin{align}
	\begin{bmatrix}
	\nabla f(\bx)\\
	\pval\bxi
	\end{bmatrix} - \begin{bmatrix}
	\nabla c(\bx)\\
	\pval  \bI
	\end{bmatrix}  (\bdual+\bz) - \begin{bmatrix}
	\nabla g(\bx)\\
	\bO
	\end{bmatrix}\bbar &=\bO\label{eqn:KKT1_Subst_ALM_unreduced1}\\
	c(\bx) + \pval  \bxi + \pmval  \bz &= \bO\,.\label{eqn:KKT1_Subst_ALM_unreduced2}
	\end{align}
	\label{eqn:KKT1_Subst_ALM_unreduced}%
\end{subequations}

We could use~\eqref{eqn:KKT1_Subst_ALM_unreduced} directly in order to form an ALM iteration. That iteration would consist of two alternating steps: 1) solving the optimality system~\eqref{eqn:KKT1_Subst_ALM_unreduced} together with \eqnKKTb for $(\bx,\bxi,\bz,\bbar,\cA)$ where $\bdual$ is fixed; 2) updating $\bdual \leftarrow \bdual + \bz$, being equivalent to $\bdual \leftarrow \bdual - \frac{1}{\pmval}\left(c(\bx)+\pval\bxi\right)$.

\subsection{Elimination of the Auxiliary Vector}
Instead, we propose to eliminate $\bxi = \bdual+\bz$ to obtain
\begin{subequations}
	\begin{align}
	\nabla f(\bx) - \nabla c(\bx)  (\bdual+\bz) - \nabla g(\bx) \bbar &=\bO\\
	c(\bx) + \pval  \bdual + (\pval+\pmval)  \bz &= \bO \label{eqn:AlternatedKKT_MALM_primal}\,.
	\end{align}
	\label{eqn:AlternatedKKT_MALM}%
\end{subequations}
As in ALM, we solve~\eqref{eqn:AlternatedKKT_MALM} and \eqnKKTb with an iteration of two alternating steps:
\begin{enumerate}
	\item Keep the value of $\bdual$ fixed, and solve~\eqref{eqn:AlternatedKKT_MALM} and \eqnKKTb for $(\bx,\bz,\bbar,\cA)$.
	\item Update $\bdual$ as $\bdual \leftarrow \bdual +\bz\,.$
\end{enumerate}
Analogous to ALM, the first step can be realized by minimizing an augmented Lagrangian function
for $\bx$ at fixed $\bdual$ subject to $\bx \in \cB$, whereas in the second step $\bz$ can be expressed in terms of $\bx$ from~\eqref{eqn:AlternatedKKT_MALM_primal}. Using this, the method can be expressed  in \abbalg~\ref{algo:MALM}, where 
\begin{align}
\Psi_{k+1}(\bx) := \cL(\bx,\bdual_{k},\bO) + \frac{0.5}{\pval + \pmval} \left\|c(\bx)+\pval  \bdual_{k}\right\|_2^2 \label{eqn:ALF}
\end{align}
is the augmented Lagrangian function, with $\cL(\bx,\bdual,\bO)\equiv f(\bx)-\bdual\t \cdot c(\bx)$.

\begin{algorithm}[tb]
	\caption{Modified Augmented Lagrangian Method}
	\label{algo:MALM}
	\begin{algorithmic}[1]
		\Procedure{MALM}{$f,c,g,\pval,\bx_0,\bdual_0,\tol$}
		\State $\pmval \leftarrow \pmval_0$
		\For{$k=1,2,3,\dots,k_\text{max}$}
			\State Compute $\bx_k$, and optionally $\bbar_k$, by solving
			\begin{align}
				&\min_{\bx \in \R^{n}}\quad\Psi_{k}(\bx)\quad\text{s.t. }\ g(\bx)\geq \bO\,. \label{eqn:ALF_Prob}
			\end{align}
			\State Update $\bdual_k \leftarrow \bdual_{k-1} - \frac{1}{\pval + \pmval} \left(c(\bx_k)+\pval \bdual_{k-1}\right)$
			\vspace{2mm}
			\If{$\|c(\bx_k)+\pval \bdual_k\|_\infty\leq\tol$}
				\State \Return $\bx_k,\bdual_k$ and optionally $\bbar_k$
			\Else
				\State Decrease $\pmval\leftarrow c_{\pmval} \pmval$ to promote convergence.
			\EndIf
		\EndFor
		\EndProcedure
	\end{algorithmic}
\end{algorithm}

\subsection{Practical Aspects}\label{sec:pracAspects}
Values that we have found work well in practice are $\tol=10^{-8}, c_\pmval=0.1, \pmval_0=0.1$. Care must be taken that $\Psi_k$ in~\eqref{eqn:ALF_Prob} is bounded below. To this end, practical methods impose box constraints $\bx_L \leq \bx \leq \bx_U$ \cite[eq.~3.2.2]{Lancelot}, expressible via $g(\bx)\geq \bO$, or a trust-region constraint $g(\bx)=\Delta^2 - \|\bx-\bx_{k-1}\|_2^2 \geq 0$ \cite[eq.~3.2.4]{Lancelot} with trust-region radius $\Delta>0$.

In order to minimize~\eqref{eqn:ALF_Prob}, one can use any numerical method for inequality constrained nonlinear minimization; e.g.\ an interior-point method like IPOPT \cite{IPOPT} or an active set method like SNOPT \cite{GilMS05}.

We refer to \cite[eq.~17.21]{NumOpt} for details on how the quasi-Newton direction for the quadratic penalty function can be computed in a more numerically stable fashion from a saddle-point linear equation system.

\subsection{Discussion}
\subsubsection{True Generalization of ALM}\label{sec:truegeneralization}
MALM is a true generalization of ALM, because they differ only by the parameter $\pval$. In particular, if $\pval=0$ then MALM in Algorithm~\ref{algo:MALM} is identical to ALM in \cite[Algorithm~3.1]{Conn95}. In contrast, when selecting $\pval>0$, we show below that MALM converges to critical points of~\eqref{eqn:CQPP} with the given~$\pval$.

\subsubsection{Benefit}
MALM solves the penalty function $\Phi_{\pval}$ in~\eqref{eqn:CQPP} by minimizing a sequence of penalty functions~$\Psi_k$. When does this make sense? If we select $\pmval \gg \pval$. Thereby, the penalty functions~$\Psi_k$ have better scaling and hence can often be minimized more efficiently in comparison to a single minimization of $\Phi_{\pval}$. The computational performance results in Section~\ref{sec:NumExp} verify this claim.

\section{Convergence Analysis}
The below analyses assume that all sub-problems~\eqref{eqn:ALF_Prob} are solved exactly, and that computations are performed in exact arithmetic. Throughout this subsection, MALM means the callback-function in \abbalg~\ref{algo:MALM}, wherein any black-box method can be used to solve~\eqref{eqn:ALF_Prob}.

\subsection{Global Convergence}
For the analyses, we consider a call of Algorithm~1 with instance $\cI:=(f,c,g,\pval,\bx_0,\bdual_0,\tol)$. MALM will create a sequence of iterates $\bx_k,\bdual_k$. 

\begin{lemma}[Equivalence]
	MALM on the instance $\cJ:=(\hat{f},\hat{c},\hat{g},0,\hbx_0,\bdual_0,\tol)$ from~\eqref{eqn:Subst} will generate the same iterates $(\bx_k,\bxi_k),\bdual_k$ as MALM on the instance $\cI$ in terms of $\bx_k,\bdual_k$.
\end{lemma}
\noindent
\textit{Proof:} By induction over $k$. Base: For $k=0$ the proposition holds by construction of the initial guesses. Step: Let the proposition hold for $k-1$. We now show that the proposition holds for $k$. The iterate $\hbx_k$ from $\cJ$ in line~4 necessarily satisfies $\nabla_{\hbx}\Psi_k(\hbx_k)- \nabla_{\hbx} \hat{g}(\hbx)\t\cdot\bbar_k=\bO$, which is equivalent to~\eqref{eqn:KKT1_Subst_ALM_unreduced} after elimination of $\bz$ by means of~\eqref{eqn:KKT1_Subst_ALM_unreduced2}. From the second component of~\eqref{eqn:KKT1_Subst_ALM_unreduced1} it follows that
\begin{align}
\bxi_k = \bxi(\bx_k,\blambda_{k-1}):= \frac{1}{\pval+\pmval} \big(\pmval\bdual_{k-1}-c(\bx_k)\big)\,.\label{eqn:bxi_func}
\end{align}
I.e.\ $\bxi_k$ is uniquely determined, hence $\hbx_k$ takes on the form $\hbx_k=(\bx_k,\bxi(\bx_k,\bdual_{k-1}))$ for some $\bx_k$. 

Substituting~\eqref{eqn:bxi_func} into the first component of~\eqref{eqn:KKT1_Subst_ALM_unreduced1} yields $\nabla\Psi_k(\bx_k)=\bO$, which is indeed identical to what $\bx_k$ in line~4 of $\cI$ satisfies. Thus, $\hbx_k=(\bx_k,\bxi(\bx_k,\bdual_{k-1}))$ with $\bx_k$ from $\cI$ is a valid $k$th iterate of $\cJ$. Finally, notice that $\bdual_k$ in $\cI,\cJ$ are identical because 
$$ -\frac{1}{\rho} \big(c(\bx_k)+\pval \bxi(\bx_k,\blambda_{k-1})\big) = -\frac{1}{\rho+\pval}\big(c(\bx_k)+\pval\blambda_{k-1}\big)\,.$$
\QED

In turn, MALM with $\pval=0$ is identical to ALM in \cite[Algorithm~3.1]{Conn95}. We can hence use the convergence result from \cite[Thm~4.6]{Conn95}:

\begin{theorem}[Global Convergence]\label{thm:globconv}
	Choose a bounded domain $\Omega \subset \R^n$. Let $\pval>0$, $c$ be bounded on $\Omega$, and let $f,c$ be twice continuously differentiable in $\Omega$, and $g$ affine. Suppose all iterates $\lbrace\bx_k\rbrace_{k \in \N_0}$ of MALM live in $\Omega$. If $\rho_0$ is sufficiently small then $\lbrace\bx_k\rbrace_{k \in \N_0}$ converges to a critical point of~\eqref{eqn:CQPP}.
\end{theorem}
\noindent
\textit{Proof:} \cite[Thm~4.6]{Conn95} shows convergence of ALM for $\cJ$ under four assumptions (AS1)-(AS4). It suffices to show that $\hat{f},\hat{c},\hat{g}$ satisfy these assumptions.

Feasibility \cite[AS1]{Conn95} of~\eqref{eqn:Subst} holds naturally by $\bxi= \frac{-1}{\pval} \,c(\bx)$.
Twice continuous differentiability \cite[AS2]{Conn95} of $\hat{f},\hat{c}$ holds per requirement.
Boundedness \cite[AS3]{Conn95} of all $\hbx_k \in \Omega \times c(\Omega)$ follows from boundedness of $\Omega$ and $c$ on $\Omega$.

The last assumption \cite[AS4]{Conn95} is more technical. Since $g$ is affine, we can express $g(\bx)=\bA\cdot\bx-\bb$, and likewise $\hat{g}(\hbx)=\hbA\cdot\hbx-\bb$, where $\hbA = [\bA\ \bO]$. We define the matrix $\bZ$ of orthonormal columns that span the null-space of $\hbA_\cA$, i.e.\ the matrix of sub-rows of $\hbA$ of the active constraints at $\hbx$. (AS4) requires $\nabla\hat{c}(\hbx)\cdot\bZ$ to be of column rank $\geq m$. Due to the special structure of $\hbA$, we see that $\bZ$ has a structure like
\begin{align*}
\bZ = \begin{bmatrix}
\begin{matrix}
\bO\\
\bI
\end{matrix}&
\begin{matrix}
\dots\\
\dots
\end{matrix}
\end{bmatrix}\,.
\end{align*}
Since $\nabla_{\hbx}\hat{c}(\hbx)\t=[\nabla c(\hbx)\t \ \pval\bI]$ has full row rank, the rank of $\nabla\hat{c}(\hbx)\t\cdot\bZ$ is bounded below by the number of columns of $\bZ$, i.e.\ bounded below by $m$. \QED

Some of the requirements in Theorem~\ref{thm:globconv} may be forcible\todo{3d}: Section~\ref{sec:pracAspects} explains how $\cB$ can be bounded. In this case, choosing $\Omega=\cB$ yields $\lbrace \bx_k\rbrace \subset \Omega$. Also, $c$ may be bounded over $\Omega$ by approximating $c(\bx)$ with $\arctan\big(c(\bx)\big)$. If $\|c(\bx)\|_2$ is very small at the minimizer of~\eqref{eqn:CQPP} then the approximation error of $\arctan$ is negligible. \todo{2p}To make $g$ affine, several practical ALM implementations (Lancelot, MINOS) convert inequalities to equalities via the addition of slack variables $\bs\geq 0$~\cite[Sec.~17.4]{NumOpt}. The constraints $g(\bx)-\bs=\bO$ (as in~\cite[eqn~17.47]{NumOpt}) can be merged into $c$ and scaled such that they hold tightly. Also, interior-point methods like IPOPT~\cite{IPOPT} use slacks to ensure iterates are strictly interior.

\subsection{Local Convergence}\label{app:Local}
\cite[Thm~17.6]{NumOpt} asserts linear convergence of ALM when $\nabla_{\bx\bx}^2\cL,\nabla c,\nabla_\bx\cL$ are local Lipschitz-continuous and $\pmval=const$. Section~\ref{sec:NumExp:Circle:Rate} and Figure~\ref{fig:convrate} show this. Likewise, MALM attains a linear rate in the limit when $\pmval=const$.\todo{1d} Upper bounds for these rates can be computed. In this section we prove that the rate of MALM is strictly smaller than that of ALM.

For the following result, we compare the iteration of ALM and MALM from the same initial guess $\bx_0,\bdual_0$ and the same problem-defining functions $f,c,g$. We assume that $\lbrace \bx_k \rbrace \subset \cU$, where $\cU\subset \R^n$ is an open neighborhood which contains unique local minimizers of both~\eqref{eqn:CP} and~\eqref{eqn:CQPP}.

\begin{theorem}[Local Convergence]\label{thm:locconv}
	Let $\nabla f$, $\nabla c$, $\nabla_{\bx\bx}^2 \cL$ be Lipschitz-continuous $\forall \bx \in \cU$ and let all iterates of ALM and MALM remain in $\cU$. Let the local minimizers satisfy strict complementarity. Apply ALM and MALM with fixed penalty parameter $\pmval$ to solve either problem, each starting from $\bx_k$. If both methods converge and if $\bx_k$ is sufficiently close to the local minimizer of~\eqref{eqn:CQPP}, then the linear rates of convergence of MALM and ALM satisfy the relation $C_\text{MALM}=\frac{\pmval}{\pmval+\pval} \cdot C_\text{ALM}<C_\text{ALM}$.
\end{theorem}
\noindent
\textit{Proof:}
We use the Taylor series
\begin{align*}
{}&{}\nabla \Psi_k(\bx_k,\bdual_{k-1})\\
={}&{} \bH\bx_k+\bg-\frac{1}{\pval+\pmval}\bJ\t\left(\pmval\bdual_{k-1}+\bc\right) + R_L(\bx_k,\bdual_{k-1})
\end{align*}
with $\bJ\t := \nabla c(\bx_\infty)$, $\bH := \nabla^2_{\bx\bx}\cL(\bx_\infty,\bdual_\infty,\bO) + \frac{1}{\pval+\rho}\bJ\t\bJ$, $\bc := \bJ \bx_\infty - c(\bx_\infty)$ and $\bg:=\nabla f(\bx_\infty)$ has the Lagrange remainder $\|R_L(\bx_k,\bdual_k)\|_2 \leq \frac{L}{\rho+\pval} (\|\bx_k-\bx_\infty\|_2 + \|\bdual_{k-1}-\bdual_\infty\|_2)^2$, where $L$ is the Lipschitz constant.

We now first consider the case where $p=0$, i.e.\ when there are no inequality constraints. Since $\bx_k$ is convergent by requirement, $\bH$ must be positive semi-definite and, if $\bx_\infty$ is locally unique, $\bH$ must be positive definite. Clearly, local convergence to a unique point depends quantitatively on uniqueness, hence we imply $\lambda_\text{min}(\bH)\geq \mu>0$. For the induced 2-norm it follows that $\|\bH\inv\|_2\leq \mu$, hence
\begin{align*}
&{}\left\|\bx_k-\bH\inv\left(\bg-\frac{1}{\pval+\pmval}\bJ\t(\pmval\bdual_{k-1}+\bc)\right)\right\|_2\\
\leq{}&{}\frac{L}{\mu(\rho+\pval)}\|\bdual_{k-1}-\bdual_\infty\|^2_2.
\end{align*}
Inserting the estimate for $\bx_k$ into line~5 in \abbalg~\ref{algo:MALM} gives a formula for $\bdual_k$ that only depends on $\bdual_{k-1}$:
\begin{align}
\bdual_k = \bM \cdot \bdual_{k-1} + \mathbf{f} + R_{\bdual}(\bdual_k) \label{eqn:Banach}
\end{align}
with $\bM \in \R^{m \times m}$  below, some $\mathbf{f} \in \R^m$, and $\|R_{\bdual}(\bdual_k)\|_2 \leq \frac{1}{\mu}\left(\frac{L}{\rho+\pval}\right)^2 \|\bdual_k-\bdual_\infty\|_2^2$. Rearranging  reveals
\begin{align*}
\bM &= \frac{\pmval}{\pval+\pmval}\left(\bI - \frac{1}{\pval+\pmval}\bJ\bH\inv\bJ\t\right)\,.
\end{align*}
Since \abbthm~\ref{thm:globconv} asserts convergence of $\bdual_k$, the second order terms become negligible compared to the first-order terms and can hence be ignored in the limit. Then,~\eqref{eqn:Banach} is a Banach iteration. Thus, in the limit, the rate of convergence for $\bdual_k$ is linear with contraction $\|\bM\|_2<1$. The analysis holds regardless of whether $\pval=0$ or $>0$.

We see that in the limit MALM converges faster than ALM because $\frac{\pmval}{\pval+\pmval}<1$ when $\pval>0$, whereas $\frac{\pmval}{\pval+\pmval}=1$ when $\pval=0$. Hence, in the limit $k\rightarrow \infty$ MALM yields a stronger contraction for the errors per iteration than ALM. This is in particular an advantage in cases where ALM would converge slowly. For instance, choosing $\rho=10\pval$ guarantees convergence in the limit with at least a rate of contraction of $\frac{\pmval}{\pval+\pmval}<0.91$\,.

\balance

From the above, when dropping the Lagrange remainder terms, we can 
identify the local rate of convergence by that of the following quadratic model iteration: 1) Solve
\begin{align*}
\operatornamewithlimits{min}_{x} \frac{1}{2}\bx\t\bH\bx + \left(\bg+\frac{1}{\pval+\pmval}\bJ\t(\bJ\bx_{k-1}-\bc-\pval\bdual_{k-1})\right)\t\bx\,.
\end{align*}
2) Update $\bdual_k:=\bdual_{k-1}-\frac{1}{\pval+\pmval}(\bJ\bx_k-\bc-\pval\bdual_{k-1})$.

We  discuss the case when $p>0$, i.e.\ when inequality constraints are present. We use our assumption on strict complementarity, i.e.\ $i \in \cA \Leftrightarrow \bbar_i>\beta$ for some real $\beta>0$. Since $\bx_k$ converges by requirement, $\bdual_{k-1}$ converges and thus also $\nabla \Psi_k(\bx_k)$ converges. Hence, $\bbar_k$ must converge in order to yield $\nabla_\bx \Psi_k(\bx_k)-\nabla g(\bx)\,\bbar_k=\bO$. Once $\bbar_k$ changes less than $\beta$ at some finite $k_0 \in \N$, the active set $\cA_k$ will remain unchanged $\cA_\infty$ for all subsequent iterations $k\geq k_0$. We use $\bg_\infty$ for only the active constraints of $g$ and define $\bA_\infty := \nabla \bg_\infty(\bx_\infty)\t$, $\bb_\infty:=\nabla \bg_\infty(\bx_\infty)\t \cdot \bx_\infty$; hence $g_\infty(\bx) = \bA_\infty \cdot \bx - \bb + \cO(\|\bx-\bx_\infty\|_2^2)$.

Given the above intermezzo, the appropriate model iteration in the limit becomes obvious: 1) Solve
\begin{align*}
&\operatornamewithlimits{min}_{x} \frac{1}{2}\bx\t\bH\bx + \left(\bg+\frac{1}{\pval+\pmval}\bJ\t(\bJ\bx_{k-1}-\bc-\pval\bdual_{k-1})\right)\t\bx\\
&\text{s.t. }\bA_\infty \cdot \bx = \bb_\infty\,.
\end{align*}
2) Update $\bdual_k:=\bdual_{k-1}-\frac{1}{\pval+\pmval}(\bJ\bx_k-\bc-\pval\bdual_{k-1})$.

This is just a projection of the iteration above. Thus, we can project the iteration for $\bx_k$ onto the nullspace of $\bA_\infty$, identifying $\bx_k = \bx_r + \bN \tbx_k$ $\forall k \geq k_0$, where $\bx_r \in \R^n$ has active set $\cA_\infty$, $\tbx_k \in \R^{n-\dim(\cA_\infty)}$ and $\bN$ is a matrix of orthogonal columns that span the nullspace of $\nabla g_\infty(\bx_\infty)$. Defining $\tbH:=\bN\t\bH\bN$, $\tbJ:=\bJ\bN$, and $\tbg,\tbc$ appropriately, we arrive at the former unconstrained quadratic model iteration form, but with $\bH,\bg,\bJ,\bc,\bx_k$ replaced by the tilded quantities. Accordingly, the Banach iteration matrix $\bM$ is replaced with the matrix 
$$ 	\tbM = \frac{\pmval}{\pval+\pmval} \left(\bI-\frac{1}{\pval+\pmval}\tbJ\tbH\inv\tbJ\t\right)\,. 	$$
The resulting contraction matrix $\tbM$ for the Banach iteration of the inequality constrained case has  a factor $\frac{\rho}{\pval+\pmval}$ in front, just like for the case when $p=0$. Thus, for $\pmval>0$ the method converges locally faster in the limit $k\rightarrow\infty$. \QED

\section{Numerical Experiments}\label{sec:NumExp}
For our tests we use two instances:~\eqref{eqn:Circle} and~\eqref{eqn:OCPdisc}. \todo{2m}Each instance will be considered once as~\eqref{eqn:CQPP} and once as~\eqref{eqn:CP}. \todo{2l}Both instances are parametric: The inconsistency of~\eqref{eqn:Circle} grows in the order of $\varepsilon$ and inconsistency of~\eqref{eqn:OCPdisc} grows in the order of the mesh size $h$. The sub-problems in~\eqref{eqn:ALF_Prob} are solved with IPOPT version 12.0.3. For tests on examples with equality constraints only, we refer to \cite{CDC2020_MALM}.

\subsection{Circle Problem}
\subsubsection{Setting}
\begin{figure}
	\centering
	\includegraphics[width=0.5\linewidth]{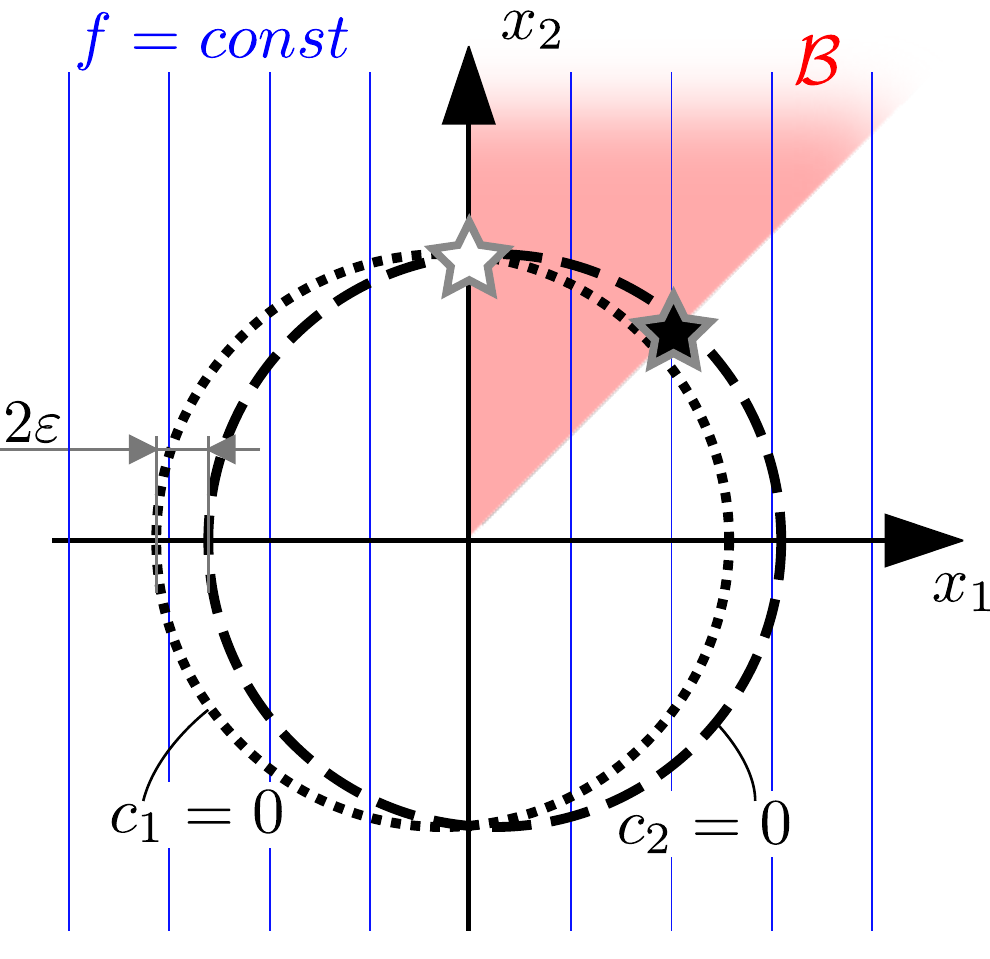}
	\caption{Geometry of the Circle Problem, with level sets of $f,c_1,c_2$ in blue solid, black dotted, and black dashed lines, respectively. The domain $\cB$ is highlighted in red. The points $\bx_A,\bx_B$ are marked as white and black star, respectively.}
	\label{fig:example1geometry}
\end{figure}

\paragraph{Initial Guess and Numerical Methods}
We use the initial guess $\bx_0=[2\,\ 1]\t$, $\bdual_0=\bO$. Fig.~\ref{fig:example1geometry} shows the instance's geometry. The figure also shows two points $\bx_A := [0\,\ \sqrt{2}]\t,\ \bx_B := [1\,\ 1]\t$.

\paragraph{Expected Minimizers}
When considering the instance as~\eqref{eqn:CP} then we expect that $\bx_A$ would be the solution. To see this, notice that $c(\bx)=\bO$ is only satisfied at $\bx=\bx_A$. When $\varepsilon \rightarrow 0$, \eqnKKT becomes ill-conditioned for $\bx_A$. Once $\varepsilon=0$, the minimizer is suddenly $\bx_B$.

When considering the instance as~\eqref{eqn:CQPP} then a point close to $\bx_B$ should be the solution unless $\varepsilon$ becomes large relative in comparison to $\pval$. To see this, note that $\bx_B$ minimizes $f$ among all points in $\cB$ that yield $\|c(\bx)\|_2$ small relative to $\pval$.

\paragraph{Scope}
Both ways~\eqref{eqn:CP} and~\eqref{eqn:CQPP} of interpreting the instance~\eqref{eqn:Circle} and both solutions $\bx_A,\bx_B$ make sense in their own right. We want to find out which solver works best for solving a respective combination $\pval,\varepsilon$.

\subsubsection{Computational Results}
We observe that all iterates of all methods remain in $\Omega=\cB \cap \lbrace \bx \in \R^2 \vert x_2 \leq 2 \rbrace$. Hence, Theorem~\ref{thm:globconv} asserts\todo{3b} that MALM and ALM converge because $f,c$ are twice continuously differentiable on $\Omega$ and $g$ is affine.

We solve the instance with MALM and QPM, for various values of $\varepsilon,\pval$, including~0. We implement QPM by solving~\eqref{eqn:CQPP} directly in IPOPT with objective $\Phi_{\pval}$. Recall that MALM=ALM for $\pval=0$ and that QPM is not applicable (n.a.) when $\pval=0$, since $\Phi_{\pval}$ is undefined.

\paragraph{Confirmation of Expected Minimizers}
We first analyze the limit points $\bx_\infty$ (which are identical for both tested methods throughout all tests) for each $\varepsilon,\pval$, by measuring the quantities
\begin{align*}
	e_A := \|\bx_\infty - \bx_A\|_2\,,\qquad e_B:=\|\bx_\infty-\bx_B\|_2\,.
\end{align*}

\abbTab~\ref{tab:Exp_Circ_conv} shows the quantities $e_A,e_B$ for respective $\varepsilon,\pval$. Dividing the table into a lower left and an upper right triangle, we see that indeed solutions in the lower triangle are close to $\bx_A$ and those on the diagonal and in the upper right are close to $\bx_B$. This confirms that solutions of~\eqref{eqn:CP} and~\eqref{eqn:CQPP} can be very distinct and the latter depend on the value of $\pval$.

\begin{table}
	\caption{Solution of the Circle Problem with respect to $\varepsilon,\pval$. Smaller values mean closer convergence to either point. Cells in the lower left converge to $\bx_A$, cells in the upper right to $\bx_B$.}
	\label{tab:Exp_Circ_conv}
	\centering
	\begin{tabular}{cl||c|c|c|c|c|}    \cline{3-7}\multicolumn{1}{l}{}&&\multicolumn{5}{c|}{$\varepsilon$}\\ \cline{2-7}  
		\multicolumn{1}{l|}{}& $\begin{matrix}e_A\\e_B\end{matrix}$			&$1.0\text{e--}1$											&$1.0\text{e--}2$ 				 								&$1.0\text{e--}4$ 												&$1.0\text{e--}6$ 												&$0.0$\\ \hline\hline
		\multicolumn{1}{|c|}{\multirow{5}{*}{$\pval$}} 	& $1.0\text{e--}1$	&$\begin{matrix}1.1\text{e--}0\\2.6\text{e--}3\end{matrix}$	&$\begin{matrix}1.1\text{e+}0	\\4.3\text{e--}3\end{matrix}$	&$\begin{matrix}1.1\text{e+}0	\\4.4\text{e--}3\end{matrix}$	&$\begin{matrix}1.1\text{e+}0	\\4.4\text{e--}3\end{matrix}$ 	&$\begin{matrix}1.1\text{e+}0	\\4.4\text{e--}3\end{matrix}$\\ \cline{2-7}
		\multicolumn{1}{|c|}{\multirow{5}{*}{}}  		& $1.0\text{e--}2$	&$\begin{matrix}1.2\text{e--}1\\9.6\text{e--}1\end{matrix}$	&$\begin{matrix}1.1\text{e+}0	\\3.7\text{e--}4\end{matrix}$	&$\begin{matrix}1.1\text{e+}0	\\4.4\text{e--}4\end{matrix}$	&$\begin{matrix}1.1\text{e+}0	\\4.4\text{e--}4\end{matrix}$ 	&$\begin{matrix}1.1\text{e+}0	\\4.4\text{e--}4\end{matrix}$\\ \cline{2-7}
		\multicolumn{1}{|c|}{\multirow{5}{*}{}} 		& $1.0\text{e--}4$	&$\begin{matrix}3.8\text{e--}3\\1.1\text{e+}0\end{matrix}$	&$\begin{matrix}1.2\text{e--}1	\\9.6\text{e--}1\end{matrix}$	&$\begin{matrix}1.1\text{e+}0	\\4.4\text{e--}6\end{matrix}$	&$\begin{matrix}1.1\text{e+}0	\\4.4\text{e--}6\end{matrix}$ 	&$\begin{matrix}1.1\text{e+}0	\\4.4\text{e--}6\end{matrix}$\\ \cline{2-7}
		\multicolumn{1}{|c|}{\multirow{5}{*}{}} 		& $1.0\text{e--}6$	&$\begin{matrix}3.5\text{e--}3\\1.1\text{e+}0\end{matrix}$	&$\begin{matrix}1.3\text{e--}3	\\1.1\text{e+}0\end{matrix}$	&$\begin{matrix}1.1\text{e+}0	\\3.7\text{e--}8\end{matrix}$	&$\begin{matrix}1.1\text{e+}0	\\4.4\text{e--}8\end{matrix}$ 	&$\begin{matrix}1.1\text{e+}0	\\4.4\text{e--}8\end{matrix}$\\ \cline{2-7}
		\multicolumn{1}{|c|}{\multirow{5}{*}{}} 	    & $1.0\text{e--}8$	&$\begin{matrix}3.5\text{e--}3\\1.1\text{e+}0\end{matrix}$	&$\begin{matrix}3.7\text{e--}5	\\1.1\text{e+}0\end{matrix}$	&$\begin{matrix}1.2\text{e--}1	\\9.6\text{e--}1\end{matrix}$	&$\begin{matrix}1.3\text{+}0	\\7.1\text{e--}9\end{matrix}$ 	&$\begin{matrix}1.3\text{e+}0	\\7.1\text{e--}9\end{matrix}$\\ \cline{2-7}
		\multicolumn{1}{|c|}{\multirow{5}{*}{}} 	    & $0.0$	&$\begin{matrix}3.5\text{e--}3\\1.1\text{e+}0\end{matrix}$	&$\begin{matrix}3.7\text{e--}5	\\1.1\text{e+}0\end{matrix}$	&$\begin{matrix}1.2\text{e--}1	\\9.6\text{e--}1\end{matrix}$	&$\begin{matrix}1.3\text{+}0	\\7.1\text{e--}9\end{matrix}$ 	&$\begin{matrix}1.3\text{e+}0	\\0.0\end{matrix}$\\ \hline
	\end{tabular}
\end{table}

\paragraph{Computational Performance}
\abbTab~\ref{tab:Exp_Circ_iter} shows the sum of the number of all inner iterations of QPM and MALM for respective $\varepsilon,\pval$. We see a trend for each of the two methods: QPM converges in a few iterations when $\pval$ is moderate. However, when $\varepsilon,\pval$ both decrease, the iteration count blows up. The trend for MALM is different. MALM converges reliably for all $\varepsilon,\pval$ in the upper right triangle, including those where $\varepsilon,\pval$ are very small.

The last row of \abbTab~\ref{tab:Exp_Circ_iter} shows ALM. ALM converges quickly to $\bx_B$ when $\varepsilon=0$. In contrast, when $\varepsilon\neq 0$ then ALM should converge to $\bx_A$ but its iteration count blows up for small $\varepsilon>0$. In two instances ALM did not converge (n.c.) within $1000$ iterations. In conclusion, ALM is inefficient when $c$ has small inconsistencies.

\begin{table}
	\caption{Total number of IPOPT iterations for MALM and QPM for the Circle Problem with respect to $\varepsilon,\pval$.  Fewer iterations mean better computational efficiency; highlighting best in slanted (QPM) or bold (MALM).}
	\label{tab:Exp_Circ_iter}
	\centering
	\begin{tabular}{cl||c|c|c|c|c|}    \cline{3-7}\multicolumn{1}{l}{}&&\multicolumn{5}{c|}{$\varepsilon$}\\ \cline{2-7}  
		\multicolumn{1}{l|}{}&$\begin{matrix}\#_\text{MALM}\\\#_\text{QPM}\end{matrix}$& $1.0\text{e--}1$& $1.0\text{e--}2$& $1.0\text{e--}4$& $1.0\text{e--}6$& $0.0$\\ \hline\hline
		\multicolumn{1}{|c|}{\multirow{5}{*}{$\pval$}} & $1.0\text{e--}1$	& $\begin{matrix}\text{28}\\ \textsl{14}\end{matrix}$ 	& $\begin{matrix}\text{22}\\ \textsl{13}\end{matrix}$& 	$\begin{matrix}\text{22}\\ 		\textsl{13}\end{matrix}$&  	$\begin{matrix}\text{19}\\ \textsl{13}\end{matrix}$& 		$\begin{matrix}\text{19}\\ \textsl{13}\end{matrix}$\\ \cline{2-7}
		\multicolumn{1}{|c|}{\multirow{5}{*}{}} & $1.0\text{e--}2$ 			& $\begin{matrix}\text{36}\\ \textsl{12}\end{matrix}$ 	& $\begin{matrix}\text{28}\\ \textsl{16}\end{matrix}$& 	$\begin{matrix}\text{16}\\ 		\text{16}\end{matrix}$& 	$\begin{matrix}\text{23}\\ \textsl{16}\end{matrix}$& 		$\begin{matrix}\text{20}\\ \textsl{16}\end{matrix}$\\ \cline{2-7}
		\multicolumn{1}{|c|}{\multirow{5}{*}{}} & $1.0\text{e--}4$ 			& $\begin{matrix}\text{21}\\ \textsl{16}\end{matrix}$ 	& $\begin{matrix}\text{56}\\ \textsl{36}\end{matrix}$& 	$\begin{matrix}\textbf{32}\\ 		\text{43}\end{matrix}$& 	$\begin{matrix}\textbf{29}\\ \text{43}\end{matrix}$& 		$\begin{matrix}\textbf{23}\\ \text{43}\end{matrix}$\\ \cline{2-7}
		\multicolumn{1}{|c|}{\multirow{5}{*}{}} & $1.0\text{e--}6$ 			& $\begin{matrix}\text{29}\\ \textsl{16}\end{matrix}$ 	& $\begin{matrix}\text{68}\\ \textsl{35}\end{matrix}$& 	$\begin{matrix}\textbf{45}\\ 		\text{138}\end{matrix}$& 	$\begin{matrix}\textbf{39}\\ \text{134}\end{matrix}$& 	$\begin{matrix}\textbf{31}\\ \text{134}\end{matrix}$\\ \cline{2-7}
		\multicolumn{1}{|c|}{\multirow{5}{*}{}} & $1.0\text{e--}8$			& $\begin{matrix}\text{34}\\ \text{n.~c.}\end{matrix}$	& $\begin{matrix}\text{60}\\ \text{n.~c.}\end{matrix}$& $\begin{matrix}\text{n.~c.}\\ 	\text{n.~c.}\end{matrix}$& 	$\begin{matrix}\textbf{52}\\ \text{429}\end{matrix}$& 	$\begin{matrix}\textbf{40}\\ \text{374}\end{matrix}$\\ \cline{2-7}
		\multicolumn{1}{|c|}{\multirow{5}{*}{}} & $0.0$ 					& $\begin{matrix}\text{34}\\ \text{n.~a.}\end{matrix}$	& $\begin{matrix}\text{60}\\ \text{n.~a.}\end{matrix}$& $\begin{matrix}\text{n.~c.}\\ 	\text{n.~a.}\end{matrix}$& 	$\begin{matrix}\text{52}\\ \text{n.~a.}\end{matrix}$& 	$\begin{matrix}\text{40}\\ \text{n.~a.}\end{matrix}$\\ \hline
	\end{tabular}
\end{table}

\begin{figure}
	\centering
	\includegraphics[width=0.7\linewidth]{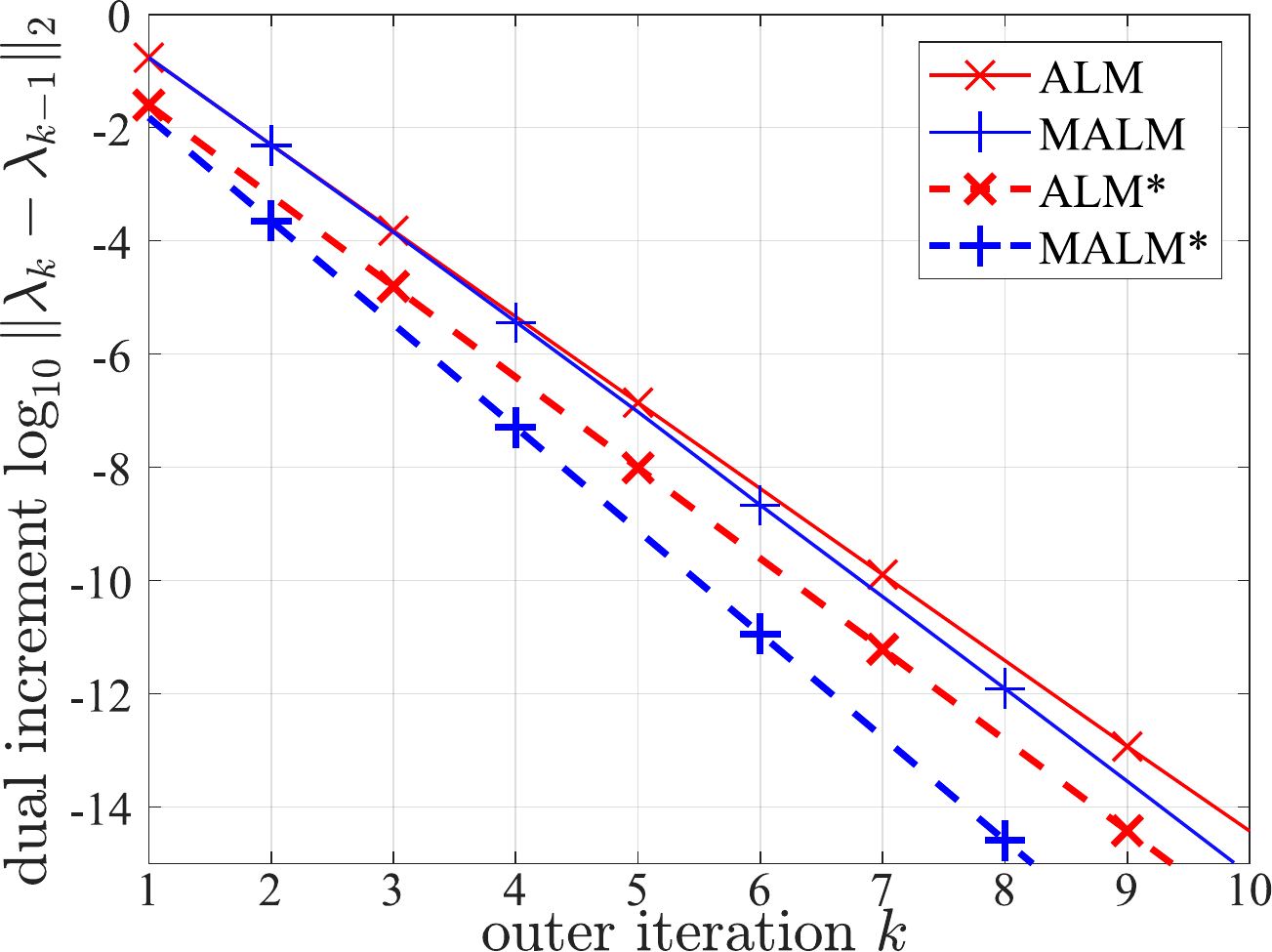}
	\caption{Comparison of convergence for ALM and MALM for the circle problem. Solid lines are measured convergence rates. Dotted lines indicate the theoretical rates of ALM and MALM in the limit $k \rightarrow \infty$.}
	\label{fig:convrate}
\end{figure}

\paragraph{Rate-of-Convergence Comparison}\label{sec:NumExp:Circle:Rate}
We compare the rate of convergence of MALM and ALM to the theoretical prediction from Theorem~\ref{thm:locconv}\todo{3a}. We use $\varepsilon=0$ and $\rho=1$\,. MALM solves~\eqref{eqn:CQPP} with $\pval=10^{-1}$ whereas ALM solves~\eqref{eqn:CP}. Both minimizers are close to $\bx_B$. Figure~\ref{fig:convrate} plots $\|\bdual_{k}-\bdual_{k-1}\|_2$ for both methods over the outer iteration index $k$ of Algorithm~\ref{algo:MALM}. We observe convergence at linear rates. We see that both methods converge in very few outer iterations to the order of machine epsilon. At $k\geq 9$ both methods have roughly attained their limit convergence rates.

\subsection{Integral Penalty-Discretization for Optimal Control}
\subsubsection{Setting}
\paragraph{Initial Guess and Solvers}
We solve the instance~\eqref{eqn:OCPdisc} with MALM and QPM for various values of $h,\pval$ from $\bx_0=\bO,\bdual_0=\bO$. Recall that $h$ is the mesh size and $\|c(\bx)\|_2^2\equiv \int_0^5 (y_h^2/2 + u_h-\dot{y}_h)^2\,\mathrm{d}t$.

For this example, using the bounded domain $\Omega=\cB$, Theorem~\ref{thm:globconv} asserts a priori\todo{3c} that MALM converges because $f,c$ are twice continuously differentiable on $\Omega$ and $g$ is affine.

\paragraph{Expected Minimizers}
We expect that the optimality gap and feasibility residual\todo{2g}
\newcommand{\symgap}{\delta J}
\begin{align*}
	\delta J= f(\bx)-J(y^\star,u^\star)\,,\quad r= \|c(\bx)\|_2^2
\end{align*}
both converge for increasing mesh sizes $N$ when choosing $\pval \in \cO(h)$; cf. discussion in Section~\ref{sec:Intro:Motiv}. For $\pval$ too large, $r$ should not converge and for $\pval$ too small $\delta J$ should not converge. To see this, notice that $\Phi_{\pval}=J+\frac{1}{2 \pval} r$; thus, minimization of $\Phi_{\pval}$ only strikes a balance between minimizing both terms when $\pval$ is chosen in the right order of magnitude

\paragraph{Scope}
We expect that again QPM will be faster than MALM when $\pval$ is moderate and vice versa when $\pval$ is very small.
We shall also try ALM (i.e., Algorithm~\ref{algo:MALM} with $\pval=0$) but just for completeness, because this will not converge to the optimal control solution.

\subsubsection{Computational Results}
\paragraph{Confirmation of Expected Minimizers}
\abbTab~\ref{tab:Exp_OCP_conv} shows the quantities $\symgap,r$ for respective $h,\pval$. Dividing the table into a lower left and an upper right triangle, we find our expected minimizers confirmed: solutions in the lower left of the table\todo{2k} achieve small $r$ but large $\delta J$, whereas solutions in the upper right of the table are not sufficiently feasible with respect to the path-constraints. For a given mesh size $h$, the most accurate control solutions are found on the diagonal cells of the table.\todo{2h}

\begin{table}
	\caption{Solution of the Optimal Control Problem with respect to $N,\pval$. For a given mesh size $N$, the value for $\pval$ is suitable when $\delta J$ (optimality gap) and $r$ (feasibility residual) have similar magnitude.}
	\label{tab:Exp_OCP_conv}
	\centering
	\begin{tabular}{cc||c|c|c|c|c|}    \cline{3-7}\multicolumn{1}{l}{}&&\multicolumn{5}{c|}{$h$}\\ \cline{2-7}  
		\multicolumn{1}{l|}{}& $\begin{matrix}\symgap\\r\end{matrix}$ 			&$1.0\text{e--}1$ 																	& $2.0\text{e--}2$ 																& $1.0\text{e--}2$ 																	& $5.0\text{e--}3$ 																& $2.5\text{e--}3$\\ \hline\hline
		\multicolumn{1}{|c|}{\multirow{5}{*}{$\pval$}} 	& $1.0\text{e--}2$		&$\begin{matrix}\text{-}1.7\text{e--}1\\9.2\text{e--}2\end{matrix}$		&$\begin{matrix}\text{-}1.7\text{e--}1	\\4.3\text{e--}2\end{matrix}$ 	&$\begin{matrix}\text{-}1.7\text{e--}1	\\4.3\text{e--}2\end{matrix}$		&$\begin{matrix}\text{-}1.7\text{e--}1\\4.2\text{e--}2\end{matrix}$		&$\begin{matrix}\text{-}1.7\text{e--}1\\4.3\text{e--}2\end{matrix}$\\ \cline{2-7}
		\multicolumn{1}{|c|}{\multirow{5}{*}{}} 		& $1.0\text{e--}3$ 		&$\begin{matrix}        4.3\text{e--}2\\2.7\text{e--}2\end{matrix}$		&$\begin{matrix}\text{-}8.6\text{e--}3	\\5.3\text{e--}3\end{matrix}$ 	&$\begin{matrix}\text{-}9.6\text{e--}3	\\4.6\text{e--}3\end{matrix}$		&$\begin{matrix}\text{-}9.8\text{e--}3\\4.4\text{e--}3\end{matrix}$		&$\begin{matrix}\text{-}9.9\text{e--}3\\4.4\text{e--}3\end{matrix}$\\ \cline{2-7}
		\multicolumn{1}{|c|}{\multirow{5}{*}{}} 		& $1.0\text{e--}4$ 		&$\begin{matrix}        7.6\text{e--}2\\8.8\text{e--}3\end{matrix}$		&$\begin{matrix}1.9\text{e--}2 			\\6.1\text{e--}4\end{matrix}$ 	&$\begin{matrix}1.2\text{e--}2			\\6.1\text{e--}4\end{matrix}$		&$\begin{matrix}\text{-}8.7\text{e--}3\\6.0\text{e--}4\end{matrix}$		&$\begin{matrix}\text{-}7.2\text{e--}3\\5.7\text{e--}4\end{matrix}$\\ \cline{2-7}
		\multicolumn{1}{|c|}{\multirow{5}{*}{}} 		& $1.0\text{e--}5$ 		&$\begin{matrix}        7.9\text{e--}2\\2.5\text{e--}3\end{matrix}$		&$\begin{matrix}2.3\text{e--}2			\\6.4\text{e--}5\end{matrix}$ 	&$\begin{matrix}1.6\text{e--}2			\\6.2\text{e--}5\end{matrix}$		&$\begin{matrix}        1.2\text{e--}2\\6.3\text{e--}5\end{matrix}$		&$\begin{matrix}        1.0\text{e--}2\\6.6\text{e--}5\end{matrix}$\\ \cline{2-7}
		\multicolumn{1}{|c|}{\multirow{5}{*}{}} 		& $1.0\text{e--}6$ 		&$\begin{matrix}        8.0\text{e--}2\\6.5\text{e--}4\end{matrix}$		&$\begin{matrix}2.3\text{e--}2			\\1.8\text{e--}5\end{matrix}$ 	&$\begin{matrix}1.6\text{e--}2			\\7.5\text{e--}6\end{matrix}$		&$\begin{matrix}        1.2\text{e--}2\\1.1\text{e--}5\end{matrix}$		&$\begin{matrix}        1.0\text{e--}2\\1.5\text{e--}5\end{matrix}$\\ \cline{2-7}
		\multicolumn{1}{|c|}{\multirow{5}{*}{}} 		& $0.0$ 				&$\begin{matrix}        7.2\text{e+}0\\0.0\end{matrix}$ 				&$\begin{matrix}7.2\text{e+}0			\\0.0\end{matrix}$ 				&$\begin{matrix}7.2\text{e+}0			\\0.0\end{matrix}$					&$\begin{matrix}        7.2\text{e+}0\\0.0           \end{matrix}$		&$\begin{matrix}        7.2\text{e+}0\\0.0           \end{matrix}$\\ \hline
	\end{tabular}
\end{table}

\paragraph{Computational Performance}
\abbTab~\ref{tab:Exp_OCP_iter} shows the sum of the number of all inner iterations of MALM and QPM for respective $h,\pval$. We see the same trend as for the circle problem: QPM converges faster than MALM when $\pval$ is moderate and vice versa when $\pval$ is small. We underline that MALM converges reliably for all $h,\pval$ in the upper right triangle, including those where $h,\pval$ are very small. Needless to say, accurate numerical optimal control solutions require $h,\pval$ very small; thus MALM seems very attractive for solving these classes of problems.

The last row shows that ALM does not convergence (n.c.) within $500$ iterations for any mesh size.

\begin{table}
	\caption{Total number of IPOPT iterations for MALM and QPM for the Optimal Control Problem with respect to $N,\pval$. Fewer iterations mean better computational efficiency; highlighting best in slanted (QPM) or bold (MALM).}
	\label{tab:Exp_OCP_iter}
	\centering
	\begin{tabular}{cc||c|c|c|c|c|}    \cline{3-7}\multicolumn{1}{l}{}&&\multicolumn{5}{c|}{$h$}\\ \cline{2-7}  
		\multicolumn{1}{l|}{}&$\begin{matrix}\#_\text{MALM}\\\#_\text{QPM}\end{matrix}$& 	$1.0\text{e--}1$ 																	& $2.0\text{e--}2$ 																& $1.0\text{e--}2$ 																	& $5.0\text{e--}3$ 																& $2.5\text{e--}3$\\ \hline\hline
		\multicolumn{1}{|c|}{\multirow{5}{*}{$\pval$}} 	& $1.0\text{e--}2$& 				$\begin{matrix}\text{39}\\ 		\textsl{17}\end{matrix}$& 	$\begin{matrix}\text{53}\\ \textsl{24}\end{matrix}$& 			$\begin{matrix}\text{64}\\    \textsl{23}\end{matrix}$& 		$\begin{matrix}\text{69}\\    \textsl{26}\end{matrix}$& 		$\begin{matrix}\text{69}\\    \textsl{44}\end{matrix}$\\ \cline{2-7}
		\multicolumn{1}{|c|}{\multirow{5}{*}{}}			& $1.0\text{e--}3$& 				$\begin{matrix}\text{47}\\ 		\textsl{38}\end{matrix}$& 	$\begin{matrix}\text{63}\\ \textsl{44}\end{matrix}$& 			$\begin{matrix}\text{74}\\    \textsl{43}\end{matrix}$& 		$\begin{matrix}\text{90}\\    \textsl{44}\end{matrix}$& 		$\begin{matrix}\text{89}\\    \textsl{68}\end{matrix}$\\ \cline{2-7}
		\multicolumn{1}{|c|}{\multirow{5}{*}{}}  		& $1.0\text{e--}4$& 				$\begin{matrix}\text{61}\\ 		\textsl{57}\end{matrix}$& 	$\begin{matrix}\textbf{66}\\ \text{92}\end{matrix}$& 			$\begin{matrix}\textbf{80}\\    \text{85}\end{matrix}$& 		$\begin{matrix}\textbf{93}\\    \text{133}\end{matrix}$& 		$\begin{matrix}\textbf{121}\\   \text{124}\end{matrix}$\\ \cline{2-7}
		\multicolumn{1}{|c|}{\multirow{5}{*}{}} 		& $1.0\text{e--}5$& 				$\begin{matrix}\textbf{64}\\ 		\text{102}\end{matrix}$& 	$\begin{matrix}\textbf{78}\\ \text{161}\end{matrix}$& 		$\begin{matrix}\textbf{80}\\    \text{266}\end{matrix}$& 		$\begin{matrix}\textbf{93}\\    \text{242}\end{matrix}$& 		$\begin{matrix}\textbf{113}\\   \text{252}\end{matrix}$\\ \cline{2-7}
		\multicolumn{1}{|c|}{\multirow{5}{*}{}} 		& $1.0\text{e--}6$& 				$\begin{matrix}\textbf{81}\\ 		\text{163}\end{matrix}$& 	$\begin{matrix}\textbf{78}\\ \text{222}\end{matrix}$& 		$\begin{matrix}\textbf{110}\\   \text{328}\end{matrix}$& 		$\begin{matrix}\textbf{97}\\    \text{278}\end{matrix}$& 		$\begin{matrix}\textbf{126}\\   \text{224}\end{matrix}$\\ \cline{2-7}
		\multicolumn{1}{|c|}{\multirow{5}{*}{}} 		& $0.0$& 							$\begin{matrix}\text{n.~c.}\\ 	\text{n.~a.}\end{matrix}$& 	$\begin{matrix}\text{n.~c.}\\ \text{n.~a.}\end{matrix}$& 	$\begin{matrix}\text{n.~c.}\\ \text{n.~a.}\end{matrix}$& 	$\begin{matrix}\text{n.~c.}\\ \text{n.~a.}\end{matrix}$& 	$\begin{matrix}\text{n.~c.}\\ \text{n.~a.}\end{matrix}$\\ \hline
	\end{tabular}
\end{table}

\section{Conclusions}
We presented a modified augmented Lagrangian method (MALM), generalized to non-convex optimization problems with additional inequality constraints. We proved global convergence for our generalized method when the inequalities are affine. A local rate-of-convergence result shows that MALM inherits all the local convergence results of ALM while the regularization in $\pval>0$ also yields a slight benefit to its rate of local convergence in the iteration limit.

Our numerical experiments demonstrate that MALM outperforms QPM when minimizing
quadratic penalty programs~\eqref{eqn:CQPP} in those situations where $\pval$ is very small, in a similar manner as ALM outperforms QPM when solving equality constrained programs~\eqref{eqn:CP}. The experiments further show that ALM cannot solve~\eqref{eqn:CQPP}, but solves~\eqref{eqn:CP} instead. Hence, MALM is the best candidate for solving~\eqref{eqn:CQPP} when $\pval$ is very small.

In this paper we have assumed that the sub-problems~\eqref{eqn:ALF_Prob} are solved to high accuracy. Future work could extend the approach to inexact iterations and sub-iterations to mild tolerances. This could reduce computations at sub-iterations where the dual is far from converged. Another open subject is the extension of global convergence analysis to the cases when $g$ is convex nonlinear or non-convex nonlinear.


\appendices


\bibliographystyle{plain} 		
\bibliography{MALMconf_refs_Martin_Eric}

\end{document}